\documentstyle{amsppt}
\pagewidth{5.05in} \pageheight{8in} \NoRunningHeads \magnification
= 1200

\topmatter
\title On Indefinite Special Lagrangian Submanifolds in Indefinite Complex Euclidean Spaces\endtitle
\author Yuxin Dong \endauthor
\thanks {Supported by Zhongdian grant of NSFC}
\endthanks
\abstract {In this paper, we show that the calibrated method can
also be used to detect indefinite minimal Lagrangian submanifolds
in $C_k^m$. We introduce the notion of indefinite special
Lagrangian submanifolds in $C_k^m$ and generalize the well-known
work of Harvey-Lawson to the indefinite case.}
\endabstract
, {53C50}, {53C80}
\endsubjclass
\endtopmatter
\document
\heading{\bf 1. Introduction}
\endheading
\vskip 0.3 true cm

In their celebrated paper [HL], Harvey and Lawson introduced four
types of calibrated geometries, which have been of growing
interest over the past ten years. Calibrated submanifolds are
distinguished classes of minimal submanifolds, which are
volume-minimizing in their respective homology classes. Special
Lagrangian submanifolds are one type of the calibrated
submanifolds, which may be defined in $C^m$ or a Calabi-Yau
$m-$fold. Due to their importance in Mirror symmetry, special
Lagrangian submanifolds have received much attentions in recent
years (cf. [Jo1] and the references therein).

On the other hand, the theory of classical strings tells us that,
during the time evolution, a string sweeps out a timelike minimal
surface $\Sigma $ in a space-time. There have already been many
works on timelike minimal surfaces (cf. [Gu1,2,3], [De1,2], [Mi]
and [IT]), and some works on timelike minimal submanifolds as well
(cf. [AAI], [Bre] and [Li]). Timelike minimal submanifolds may be
viewed as simple but nontrivial examples of $D-$branes, which play
an important role in string theory too. More general, we may
investigate so-called indefinite minimal submanifolds. Notice that
these submanifolds, including timelike minimal submanifolds, are
relatively unstudied in contrast to minimal submanifolds in
Euclidean spaces. It would be interesting to explore the relation
and difference between indefinite minimal submanifolds and
classical minimal submanifolds. We also hope to find more
nontrivial examples of indefinite minimal submanifolds.

In this paper, we show that the special Lagrangian calibration on
$C^m$ introduced in [HL] can also be used to detect indefinite
minimal Lagrangian submanifolds in an indefinite complex Euclidean
space $C_k^m$. Recall that Harvey and Lawson [HL] used two methods
to show that a special Lagrangian submanifold in $C^m$ is minimal:
one is the volume comparison argument, the other is to
differentiate the calibration along the submanifold. We observe
that the second method is still valid in the indefinite case (see
Proposition 2.1). Then we may introduce the notion of indefinite
special Lagrangian submanifolds in $C_k^m$ and generalize most
results on special Lagrangian submanifolds in [HL] to the
indefinite case. Notably the potential function of a graphic
indefinite special Lagrangian submanifold satisfies a `hyperbolic
equation'. When $k=1$ and $m\geq 3$, the equation becomes a fully
nonlinear hyperbolic equation. We shall discuss both local and
global Cauchy problems for this nonlinear hyperbolic equation. As
in [HL], the methods of the moment map and the normal bundle
construction will also be used to construct some explicit
indefinite special Lagrangian submanifolds in $C_k^m$. It turns
out that indefinite special Lagrangian submanifolds exist in
abundance. In [Hi], Hitchin introduced another kind of special
Lagrangian submanifolds in $(C^m=R^m\times R^m,dxdy)$, which are
actually spacelike with respect to the metric $dxdy$. Jost and Xin
[JX] then showed that such a submanifold has mean curvature
$H\equiv 0$. Later Warren [Wa] gave a spacelike calibrated
characterization to show that a special Lagrangian graph in
Hitchin's sense has volume maximizing property. Notice that the
indefinite metric of $C_k^m$ is compatible with the complex
structure $J$ of $C^m$ in Kaherian sense while the metric $dxdy$
isn't. We may show that an indefinite special Lagrangian
submanifold in $C_k^m$ is neither volume minimizing nor volume
maximizing (see the Appendix for a more general result). Therefore
we can't use the volume comparison method in this case. The use of
the terminology `calibration' in this paper is only to emphasize
that these submanifolds are also characterized by a special closed
differential form.

\heading{\bf 2. Preliminaries}
\endheading
\vskip 0.3 true cm

Let $R_n^N$ denote the $N-$dimensional Euclidean space $R^N$
endowed with the following pseudo-Euclidean metric
$$
g_{(n,N)}=-\sum_{j=1}^ndx_j^2+\sum_{j=n+1}^Ndx_j^2\tag{1}
$$
We will call $R_n^N$ the pseudo-Euclidean space with index $n$.
Let $M$ be an indefinite submanifold in $R_n^N$, by which we mean
a submanifold whose induced metric from $R_n^N$ is non-degenerate.
The normal space $T_p^{\bot }M$ is, by definition, the orthogonal
complement of $T_pM$ in $T_pR^N$ with respect to the metric
$g_{(n,N)}$. Since the induced metric on $TM$ is non-degenerate,
we see that $T_pM\cap T_p^{\bot }M=\{0\}$ at each point of $M$ and
the induced metric on $T_p^{\bot }M$ is also non-degenerate.
Therefore we have a natural decomposition $(TR_n^N)|_M=TM\oplus
T^{\bot}M$. Let $D$, $\nabla $ denote the Levi-Civita connections
in $TR_n^N$ , $TM$ respectively and let $\nabla^{\perp}$ denote
the induced normal connection in $T^{\bot }M$. Then the formulae
of Gauss and Weingarten are given respectively by
$$
\aligned
&D_XY=\nabla _XY+h(X,Y)\\
&D_X\xi =-A_\xi X+\nabla _X^{\perp }\xi\endaligned\tag{2}
$$
for $X,Y$ tangent to $M$ and $\xi $ normal to $M$, where $h$,
$A_\xi$ are the second fundamental form and the Weingarten
transformation respectively. From (2), we easily get
$$
g_{(n,N)}(h(X,Y),\xi )=g_{(n,N)}(A_\xi X,Y)=g_{(n,N)}(X,A_\xi
Y)\tag{3}
$$
which means that $A_\xi $ is self-adjoint w.r.t. $g_{(n,N)}$.

Now we consider the complex Euclidean $m-$space $C^m$ with complex
coordinates $z_1,...,z_m$ endowed with the following
pseudo-Hermitian metric
$$
h_{(k,m)}=-\sum_{j=1}^kdz_jd\overline{z}_j+\sum_{j=k+1}^mdz_jd\overline{z}_j\tag{4}
$$
The pair $(C^m,h_{(k,m)})$ is denoted by $C_k^m$ which is called
the indefinite complex Euclidean $m-$space with index $k$.

The group of matrices in $GL(m,C)$ which leave invariant
$h_{(k,m)}$ is denoted by $U(k,m-k)$. Set
$$
I_{k,m}=\pmatrix -I_k & 0 \\
0 & I_{m-k}\endpmatrix\tag{5}
$$
Then
$$
U(k,m-k)=\{A\in GL(m,C):\overline{A}^tI_{k,m}A=I_{k,m}\}\tag{6}
$$
Write
$$
h_{k,m}=g-i\omega _{(k,m)}\tag{7}
$$
It is easy to see that $g$ is a pseudo-Euclidean metric on
$R^{2m}$ with index $2k$, which will be denoted by $g_{(2k,2m)}$.
Obviously, we have
$$
\aligned g_{(2k,2m)}(JX,JY)&=g_{(2k,2m)}(X,Y)\\
\omega_{(k,m)}(X,Y)&=g_{(2k,2m)}(JX,Y)\endaligned\tag{8}
$$
Using the canonical coordinates $(x_1,....,x_m,y_1,...,y_m)$ with
$z_j=x_j+iy_j$ ($j=1,...,m$), we may express $g_{(2k,2m)}$ and
$\omega _{(k,m)}$ as
$$
\aligned
g_{(2k,2m)}&=-\sum_{j=1}^k(dx_j^2+dy_j^2)+\sum_{j=k+1}^m(dx_j^2+dy_j^2)\\
\omega _{(k,m)}&=-\sum_{j=1}^kdx_j\wedge
dy_j+\sum_{j=k+1}^mdx_j\wedge dy_j\endaligned\tag{9}
$$
In this paper, we will investigate Lagrangian submanifolds with
respect to the symplectic form $\omega _{(k,m)}$. An
$n-$dimensional submanifold $i:M^n\hookrightarrow C_k^n$ is called
Lagrangian if $i^{*}\omega _{(k,m)}\equiv 0$ and
$i^{*}g_{(2k,2m)}$ is non-degenerate, which are equivalent to the
property that $J$ interchanges the tangent and the normal space of
$M$. Here the normal space is determined by the pseudo-Euclidean
metric $g_{(2k,2m)}$. In particular, a real $m-$plane $\zeta$ in
$C^m_k$ is Lagrangian if and only if $\omega _{(k,m)}|_\zeta =0$
and $g_{(2k,2m)}|_\zeta$ is non-degenerate. Obviously the induced
metric on a Lagrangian submanifold of an indefinite complex space
of index $k$ is a non-degenerate metric with real index $k$.

Let $M$ be a Lagrangian submanifold in $C_k^m$. Then we have
$$
\aligned
\nabla_X^{\perp }JY&=J\nabla_XY\\
A_{JY}X&=-Jh(X,Y)=A_{JX}Y\\
<h(X,Y),JZ>&=<h(Y,Z),JX>=<h(Z,X),JY>\endaligned\tag{10}
$$
for $X,Y,Z$ tangent to $M$. The mean curvature vector of $M$ is
defined by
$$
H=-\sum_{j=1}^kh(e_j,e_j)+\sum_{j=k+1}^mh(e_j,e_j)\tag{11}
$$
where $\left\{ e_j\right\}_{=1}^m$ is a Lorentz basis of $T_pM$
with $<e_j,e_l>=\varepsilon _j\delta _{jl}$
$$
\varepsilon _j=\cases -1,&\text{if $j\leq k$}\\
1,&\text{if $j\geq k+1$}\endcases\tag{12}
$$
If $H\equiv 0$, then $M$ is called minimal.
\definition{Definition 2.1} An oriented $m-$plane $\varsigma$
in $C_k^m$ is called special Lagrangian if
\newline (1) $\varsigma $ is Lagrangian w.r.t. $\omega _{(k,m)}$;
\newline (2) $\varsigma =A\varsigma _0$, where $A\in SU(k,m)$, $\varsigma
_0=R_k^m$.
\enddefinition

Define a family of holomorphic $m-$form $dz_\theta $ with $\theta
\in R$ as follows:
$$dz_\theta=e^{i\theta }dz$$
where $dz=dz_1\wedge \cdots \wedge dz_m$.

\proclaim{Proposition 2.1}A connected Lagrangian submanifold $M$
of $C_k^m$ is minimal if and only if $dz(TM)$ is a constant
complex number with norm $1$, or equivalently, $dz_\theta
(TM)\equiv 1$ for some constant phase $\theta $.
\endproclaim
\demo{Proof} Choose a local Lorentz frame field $\{e_j\}_{j=1}^m$
for $M$ such that $(\nabla e_j)_p=0$. Let $\eta _j=\frac \partial
{\partial x_j},$ $j=1,...,m$, be the standard Lorentz basis of
$R^m_k$. Then there exists a matrix $A\in U(k,m-k)$ such that
$$
e_1\wedge \cdots\wedge e_m=A(\eta _1\wedge \cdots \wedge \eta _m)
$$
So
$$
dz(e_1,...,e_m)=\det A=e^{i\theta}
$$
For any tangent $X\in T_pM$, we have
$$
\aligned X(e^{i\theta})&=\sum_{j=1}^mdz(e_1,...,h(X,e_j),....e_m)\\
&=\sum_{j=1}^m\varepsilon_j<h(X,e_j),Je_j>dz(e_1,...,Je_j,....e_m)\\
&=i\sum_{j=1}^m\varepsilon_j<h(e_j,e_j),JX>dz(e_1,...,e_j,....e_m)\\
&=i<H,JX>e^{i\theta}\endaligned
$$
Therefore $H\equiv 0$ if and only if $e^{i\theta }$ is
constant.\qed
\enddemo
\definition{Definition 2.2}Let $M$ be a Lagrangian submanifold of
$C_k^m$ with respect to $\omega _{(k,m)}$. If $dz_\theta
(TM)\equiv 1$, then $M$ is called an indefinite $\theta -$special
Lagrangian submanifold. In particular, if $\theta =0$, $M$ is
called an indefinite special Lagrangian submanifold. When $k=1$,
these Lagrangian submanifolds are called timelike too.
\enddefinition

We see that an indefinite special Lagrangian submanifold is just a
Lagrangian submanifold in $C_k^m$ whose each tangent plane is
special Lagrangian in the sense of Definition 2.1. Obviously
$dz_\theta (TM)\equiv 1$ is equivalent to $Re(dz_\theta
)(TM)\equiv 1$. Although the special Lagrangian calibration
$Re(dz_\theta )$ introduced in [HL] does not have the usual volume
property with respect to the indefinite metric, Proposition 2.1
still provides us a way to find and characterize indefinite
minimal Lagrangian submanifolds.

Obviously, $R^m=\{(x_1,...,x_m):x_i\in R\}$ is a Lagrangian plane
in $C_k^m$ whose induced metric is a pseudo-Euclidean metric with
index $k$ given by
$$
g_{(k,m)}=-\sum_{j=1}^kdx_j^2+\sum_{j=k+1}^mdx_j^2\tag{13}
$$
The group of matrices in $GL(m,R)$ which preserve the metric
$g_{(m,k)}$ is
$$
O(k,m-k)=\{A\in GL(m,R):A^tI_{k,m}A=I_{k,m}\}\tag{14}
$$

Let $G(m,2m)$ denote the Grassmann manifold of oriented real
$m-$planes in $C^m=R^m\oplus R^m$ and let $Lag(k,m)$ denote the
subset consisting of Lagrangian planes with respect to $\omega
_{(k,m)}$. Obviously the $(k,m)-$unitary group $U(k,m)$ acts on
$Lag(k,m)$ transitively. The isotropy subgroup of $U(k,m)$ at the
point $\varsigma _0=R_k^m$ is $SO(k,m-k)$ which acts diagonally on
$R_k^m\oplus R_k^m$. Thus
$$
Lag(k,m)\cong U(k,m)/SO(k,m-k)\tag{15}
$$
Notice that some real $m-$planes (e.g. $R^m$) are Lagrangian with
respect to both symplectic structures $\omega _{(k,m)}$ and
$\omega $, where $\omega $ denotes the standard symplectic
structure of $C^m$. Obviously
$$
U(k,m)\cap U(m)=U(k)\times U(m-k)
$$
So
$$
\aligned &\{P\in G(m,2m):P \text{ is Lagrangian w.r.t. $\omega
_{(k,m)}$ and $\omega$}\}\\
&=\{P\in G(m,2m):P=A\cdot R^m, A\in U(k)\times U(m-k)\}\endaligned
$$
We observe that if $M_1,M_2$ are special Lagrangian submanifolds
of $C^k$ and $C^{m-k}$ respectively (in the sense of [HL]), then
$M_1\times M_2$ is an indefinite special Lagrangian submanifold of
$C_k^m$. This is a trivial example in some sense, which is of
little interest.

In general, a Lagrangian $m-$plane w.r.t. $\omega _{(k,m)}$ is not
Lagrangian w.r.t. $\omega $ and vice versa. Let's see an example.
\example{Example 2.1}We consider two symplectic structures $\omega
_{(1,2)}$ and $\omega $ on $C^2=R^2\oplus R^2$ which are given
respectively by
$$
\aligned \omega _{(1,2)} &=-dx_1\wedge dy_1+dx_2\wedge dy_2\\
\omega &=dx_1\wedge dy_1+dx_2\wedge dy_2\endaligned
$$
Set $\varsigma_1 =span_R\{a\frac \partial {\partial x_1}+\frac
\partial {\partial x_2},\frac \partial {\partial y_1}+a\frac
\partial {\partial y_2}\}$ and $\varsigma_2 =span_R\{a\frac \partial {\partial x_1}+\frac
\partial {\partial x_2},\frac \partial {\partial y_1}-a\frac
\partial {\partial y_2}\}$ with $a\neq 0, 1$. Here the condition $a\neq 1$ is only to ensure
that the induced metric on $\varsigma_1$ from $g_{(2,4)}$ is
non-degenerate. It is easy to see that $\varsigma_1$ is Lagrangian
w.r.t. $\omega _{(1,2)}$ and not Lagrangian w.r.t. $\omega $,
while $\varsigma_2$ is just the reverse.
\endexample

\proclaim{Lemma 2.2} Suppose $\varsigma \in G(m,2m)$. Then
$\varsigma $ or $-\varsigma $ is special Lagrangian in $C_k^m$ if
and only if
\newline (1) $\varsigma $ is Lagrangian w.r.t. $\omega _{(k,m)}$;
\newline (2) $\beta (\varsigma )=0$
\newline where $\beta = Im\{dz\}$.
\endproclaim
\demo{Proof}Let $A$ be any complex linear map sending $\varsigma
_0$ to $\lambda \varsigma $ with $\lambda \in R$, i.e.,
$$
A(\eta _1\wedge \cdots \wedge \eta _m)=\lambda \varsigma
$$
where $\eta _1\wedge \cdots \wedge \eta _m=\zeta _0$. Thus we have
$$
\det A=\lambda dz(\varsigma)
$$
If $\varsigma $ is Lagrangian, then we have $dz(\varsigma
)=e^{i\theta }$. It follows that
$$
Im\{\det A\}=\lambda\sin\theta
$$
Therefore $\beta(\varsigma )=0$ if and only if $dz(\varsigma )=1$
or $-1$.\qed
\enddemo

Now we present an implicit formulation of indefinite special
Lagrangian submanifolds, which will be used later.

\proclaim{Lemma 2.3}Suppose that $f_1,...,f_m$ are smooth real
valued functions on an open set $\Omega \subset C_k^m$ and suppose
that $df_1,...,df_m$ are linearly independent at each point of
$M=\{z\in \Omega :f_1(z)=\cdots =f_m(z)=0\}$. Then the submanifold
$M$ is Lagrangian with respect to $\omega _{(k,m)}$if and only if
$$
\aligned &-\sum_{l=1}^k(\frac{\partial f_i}{\partial
x_l}\frac{\partial f_j}{\partial y_l}-\frac{\partial f_i}{\partial
y_l}\frac{\partial f_j}{\partial
x_l})+\sum_{l=k+1}^m(\frac{\partial f_i}{\partial
x_l}\frac{\partial f_j}{\partial y_l}-\frac{\partial f_i}{\partial
y_l}\frac{\partial f_j}{\partial x_l})\\
&=-2i\sum_{l=1}^k(\frac{\partial f_i}{\partial
\overline{z}_l}\frac{\partial f_j}{\partial z_l}-\frac{\partial
f_i}{\partial z_l}\frac{\partial f_j}{\partial
\overline{z}_l})+2i\sum_{l=k+1}^m(\frac{\partial f_i}{\partial
\overline{z}_l}\frac{\partial f_j}{\partial z_l}-\frac{\partial
f_i}{\partial z_l}\frac{\partial f_j}{\partial
\overline{z}_l})\endaligned\tag{16}
$$
vanish on $M$ and
$$
\det (g_{(2k,2m)}(\nabla ^gf_i,\nabla ^gf_j))\neq 0\tag{17}
$$
everywhere on $M$, where
$$\nabla
^gf_i=-\sum_{l=1}^k\frac{\partial f_i}{\partial x_l}\frac \partial
{\partial x_l}+\sum_{l=k+1}^m\frac{\partial f_i}{\partial
x_l}\frac \partial {\partial x_l}-\sum_{l=1}^k\frac{\partial
f_i}{\partial y_l}\frac
\partial {\partial y_l}+\sum_{l=k+1}^m\frac{\partial f_i}{\partial
y_l}\frac \partial {\partial y_l}
$$
are the gradient vector fields of $f_i$, $i=1,...,m$, with respect
to $g_{(2k,2m)}$.
\endproclaim
\demo{Proof}The gradient vector fields $\{\nabla
^gf_i\}_{i=1,...,m}$ are obviously linearly independent, because
$df_1,...,df_m$ are linearly independent and $g_{(2k,2m)}$ is
non-degenerate. The condition (18) ensures that the induced metric
on $M$ is non-degenerate. Since $\{\nabla ^gf_i\}_{i=1}^m$ span
the normal space at each point of $M$, the submanifold $M$ is
Lagrangian with respect to $\omega _{(k,m)}$ if and only if
$$
\omega _{(k,m)}(\nabla ^gf_i,\nabla ^gf_j)=0
$$
By a direct computation, we may derive the conclusion of this
Lemma. \qed \enddemo

\proclaim{Proposition 2.4}Suppose $M=\{z\in \Omega :f_1(z)=\cdots
=f_m(z)=0\}$ is an implicitly described Lagrangian submanifold of
 $C_k^m$. Then $M$ (with the correct orientation) is an indefinite
special Lagrangian if and only if
\newline(1) $Im\{\det_C(\partial f_j/\partial \overline{z}_l)\}=0$ on $M$ for
$m$ even;
\newline (2) $Re\{\det_C(\partial f_j/\partial
\overline{z}_l)\}=0$ on $M$ for $m$ odd.
\endproclaim
\demo{Proof} Since $M$ is Lagrangian, the tangent space of $M$ is
spanned (over $R$) by
$$
\aligned J\nabla ^gf_i &=\sum_{l=1}^k\frac{\partial f_j}{\partial
y_l}\frac \partial {\partial x_l}-\sum_{l=k+1}^m\frac{\partial
f_j}{\partial y_l}\frac \partial {\partial
x_l}-\sum_{l=1}^k\frac{\partial f_j}{\partial x_l}\frac \partial
{\partial y_l}+\sum_{l=k+1}^m\frac{\partial f_j}{\partial
x_l}\frac
\partial{\partial y_l}\\
&=\sum_{l=1}^k(\frac{\partial f_j}{\partial y_l}-i\frac{\partial
f_j}{\partial x_l})\frac
\partial {\partial x_l}-\sum_{l=k+1}^m(\frac{\partial
f_j}{\partial y_l}-i\frac{\partial f_j}{\partial x_l})\frac
\partial {\partial x_l}\\
&=(-2i\frac{\partial f_j}{\partial
\overline{z}_1},...,-2i\frac{\partial f_j}{\partial
\overline{z}_k},2i\frac{\partial f_j}{\partial
\overline{z}_{k+1}},...,2i\frac{\partial f_j}{\partial
\overline{z}_m})\endaligned
$$
where we use the natural identification of $C^m$ with $R^{2m}$. So
the complex matrix $2i(\partial f_j/\partial
\overline{z}_l)I_{k,m}$ sends $\{\frac
\partial {\partial x_1},...,\frac \partial {\partial x_m}\}$ into
the above basis for the tangent space of $M$. Hence this
Proposition follows immediately from Lemma 2.2 and Lemma 2.3.\qed
\enddemo

\heading{\bf 3. Indefinite special Lagrangian graphs}
\endheading
\vskip 0.3 true cm

First, we hope to derive the differential equation describing a
graphic indefinite special Lagrangian submanifold.

\proclaim{Lemma 3.1}Suppose $\Omega \subseteq R^m$ is open and
$f:\Omega \rightarrow R^m$ is a $C^{\infty}$ mapping. Let
$M=(x,f(x))$ be the graph of $f=(f_1,...,f_m)$ in $C_k^m$
satisfying
$$
\det \{I_{k,m}+(\frac{\partial f_l}{\partial
x_i})^tI_{k,m}(\frac{\partial f_l}{\partial x_j})\}\neq 0\tag{18}
$$
everywhere in $\Omega $. Then $M$ is Lagrangian with respect to
$\omega _{(k,m)}$ if and only if the matrix $(\partial
f^i/\partial x^j)I_{k,m}$ is symmetric. In particular, if $\Omega
$ is simply connected, then $M$ is Lagrangian with respect to
$\omega _{(k,m)}$ if and only if $f=(\nabla u)I_{k,m}$ where
$\nabla u=(u_{x_1},...,u_{x_m})$ is the gradient of some potential
function $u\in C^{\infty}(\Omega )$.
\endproclaim
\demo{Proof}It is easy to see that the induced metric on $M$ is
non-degenerate if and only if (18) holds. We may replace $f$ by
its Jacobian $f_{*}$ at some fixed point. Then
$f_{*}:R^m\rightarrow R^m$ is linear and its graph is of the form
$TM=\{x+if_{*}(x):x\in R^m\}$. By definition $TM$ is Lagrangian if
and only if $Jv\perp TM$ for all $v\in TM$ with respect to
$g_{(2k,2m)}$. Suppose $v=x+if_{*}(x)$. Then $Jv=-f_{*}(x)+ix$.
Thus $TM$ is Lagrangian if and only if $-f_{*}(x)+ix$ and
$y+if_{*}(y)$ are orthogonal for all $x,y\in R^m$, i.e.,
$$
-g_{(k,m)}(f_{*}(x),y)+g_{(k,m)}(x,f_{*}(y))=0
$$
for all $x,y\in R^m$. Write $f_{*}=A$ a $m\times m$ matrix. Then
$$
(Ax)^tI_{k,m}y=x^tI_{k,m}Ay
$$
i.e.,
$$
(AI_{k,m})^t=AI_{k,m}
$$
Thus the Jacobian of the map $fI_{k,m}$ is
$f_{*}I_{k,m}=AI_{k,m}$. Since $\Omega $ is simply connected, this
is equivalent to the existence of a potential function $u:\Omega
\rightarrow R$ with $\nabla u=fI_{k,m}$, i.e., $f=(\nabla
u)I_{k,m}$. \qed
\enddemo

For $f=(\nabla u)I_{k,m}$, we easily derive the following
$$
\aligned \det \{I+(\frac{\partial f_l}{\partial
x_i})^tI_{k,m}(\frac{\partial f_l}{\partial x_j})I_{k,m}\}
&=\det\{I+Hess(u)I_{k,m}Hess(u)I_{k,m}\}\\
&=det\{(I+iHess(u)I_{k,m})(I-iHess(u)I_{k,m})\}
\endaligned
$$
which implies that the condition (18) is equivalent to
$$
det(I+iHess(u)I_{k,m})\neq 0\tag{19}
$$
everywhere. From Lemma 2.2 and Lemma 3.1, we easily derive the
following:

\proclaim{Theorem 3.2}Suppose $u\in C^{\infty}(\Omega )$ with
$\Omega \subset R^m$. Let $M=(x,f(x))$ be the graph of $f=(\nabla
u)I_{k,m}$ in $C_k^m=R_k^m\oplus R_k^m$ satisfying (19)
everywhere. Then $M$ (with the correct orientation) is special
Lagrangian if and only if
$$
Im\{\det(I+iHess(u)I_{k,m})\}=0 \tag{20}
$$
or equivalently
$$
Im\{\det(I_{k,m}+iHess(u)\}=0
$$
\endproclaim

Let's investigate some special cases of (20). First we consider
the case $m=2$ and $k=1$. By a direct computation, we see that
(20) in this case is equivalent to
$$
u_{x_1x_1}-u_{x_2x_2}=0\tag{21}
$$
which is the one dimensional wave equation. The general smooth
solution of (21) on $R^2$ may be expressed as
$$
u=F(x_1+x_2)+G(x_1-x_2)
$$
where $F,G\in C^\infty (R)$. Consequently, (19) holds for the
graph of $f=(-u_{x_1},u_{x_2})$ if and only if
$$
4F^{\prime \prime }G^{\prime \prime }+1\neq 0\tag{22}
$$
everywhere (see also the proof of Corollary 3.4). Hence we have

\proclaim{Proposition 3.3}Let $u=F(x_1+x_2)+G(x_1-x_2)$ with
$F,G\in C^\infty (R)$. If $F$ and $G$ satisfy (22), we have a
timelike special Lagrangian surface $M=(x_1,x_2,-u_{x_1},u_{x_2})$
in $C_1^2$. Conversely, every two dimensional timelike special
Lagrangian graph is obtained in this way.
\endproclaim
\remark{Remark 3.1}By choosing any functions $F,G\in C^l(R)$ with
$l\geq 3$, we may get a $C^{l-1}$ timelike special Lagrangian
surface.
\endremark

\proclaim{Corollary 3.4}Let
$i:M=(x_1,x_2,-u_{x_1},u_{x_2})\hookrightarrow C_1^2$ be a
timelike special Lagrangian graph on $R_1^2$. Then $M$ is
conformally diffeomorphic to $R_1^2$.
\endproclaim
\demo{Proof}For the immersion
$i:M=(x_1,x_2,-u_{x_1},u_{x_2})\hookrightarrow C_1^2$, we compute
the induced metric on $M$ as follows:
$$
\aligned
di(\frac \partial {\partial x_1})&=(1,0,-u_{x_1x_1},u_{x_2x_1})\\
di(\frac\partial {\partial
x_2})&=(0,1,-u_{x_1x_2},u_{x_2x_2})\endaligned\tag{23}
$$
From (23), we have
$$
\aligned <di(\frac\partial{\partial x_1}),di(\frac \partial
{\partial x_1})>&=-1-u_{x_1x_1}^2+u_{x_1x_2}^2\\
<di(\frac\partial{\partial x_2}),di(\frac \partial {\partial
x_2})>&=1-u_{x_1x_2}^2+u_{x_2x_2}^2\\
<di(\frac \partial {\partial x_1}),di(\frac \partial {\partial
x_2})>&=-u_{x_1x_2}u_{x_1x_1}+u_{x_1x_2}u_{x_2x_2}=0\endaligned
$$
Hence the induced metric is given by
$$
ds_M^2=\lambda(-dx_1^2+dx_2^2)
$$
where $\lambda=1-u_{x_1x_2}^2+u_{x_2x_2}^2$.\qed
\enddemo

\remark{Remark 3.2}It is known that there are uncountably many
conformal structures on a simply connected Lorentz surfaces
([SW1,2]). Corollary 3.4 shows that the special Lagrangian
condition imposes a strong restriction on the conformal type of
the Lorentz graph. Recall also that the conformal Bernstein
Theorem of [Mi] states that any entire timelike minimal surface in
$R_1^3$ is $C^\infty -$conformally diffeomorphic to $R_1^2$.
\endremark

Next, for the special case $m=3$ and $k=1$, (20) becomes the
following nice form:
$$
\det (Hess(u))=\square u \tag{24}
$$
where $\square u=u_{x_1x_1}-u_{x_2x_2}-u_{x_3x_3}$.

Now we investigate the linearization of the indefinite special
Lagrangian equation at any solution. Suppose we are given an
indefinite special Lagrangian graph $M=(x,f(x))$ on $R^m$, where
$f(x)=(\nabla u(x))I_{k,m}$. Then the special Lagrangian
conditions for $M$ are
$$
\cases Re\{\det (I+if_{*})\}>0\\
Im\{\det (I+if_{*})\}=0\endcases\tag{25}
$$
For any scalar function $v$ on $R^m$, we may consider the
linearized operator
$$
\aligned L_u(v):&=Im\frac d{dt}\det
\{I+iHess(u)I_{k,m}+itHess(v)I_{k,m}\}|_{t=0}\\
&=(-1)^kIm\frac d{dt}\det
\{I_{k,m}+iHess(u)+itHess(v)\}|_{t=0}\endaligned
$$
Set $A=I_{k,m}+iHess(u)$. Notice that
$$
\det (A+itHess(v))=\det A\det (I+itA^{-1}Hess(v))
$$
Therefore
$$
\frac d{dt}|_{t=0}\det (A+itHess(v))=tr\{iA^{*}Hess(v)\}\tag{26}
$$
where $A^{*}$ denote the transposed matrix of cofactor of $A$.
Thus
$$
L_u(v)=tr\{(-1)^kRe(A^{*})Hess(v)\}\tag{27}
$$
We may diagonalize $A$ at a point $x$ so that
$$
A=diag( -1+i\lambda _1,\cdots,-1+i\lambda _k,1+i\lambda
_{k+1},\cdots,1+i\lambda _m)\tag{28}
$$
The first condition of (25) becomes
$$
(-1)^k\det A>0\tag{29}
$$
From (28) we obtain
$$
Re(A^{*})=diag(\frac{-1}{1+\lambda
_1^2},\cdots,\frac{-1}{1+\lambda _k^2},\frac 1{1+\lambda
_{k+1}^2},\cdots,\frac 1{1+\lambda _m^2})\det A \tag{30}
$$
Hence we get from (27), (29) and (30) the following:

\proclaim{Theorem 3.5}The linearization of the indefinite special
Lagrangian operator at any solution $u$ of the equation (20) is a
homogeneous second order partial differential operator
$$
L_u(v)=\sum_{ij}a^{ij}(\partial ^2u)\frac{\partial ^2v}{\partial
x_i\partial x_j}
$$
where $(a_{ij}(\partial ^2u))$ is a non-degenerate symmetric
matrix with index $k$ at each point.
\endproclaim

From Theorem 3.5, we know that Eq.(20) is an ultra-hyperbolic
equation in general. When $k=1$ and $m\geq 3$, Eq.(20) becomes a
fully nonlinear hyperbolic equation.

In the remaining part of this section, we assume that $k=1$ and
$m\geq 3$. Notice that $u\equiv 0$ is a trivial solution of (20)
and the corresponding linearization there is $-\square$, where
$\square$ is just the wave operator defined by
$$
\square v=\frac{\partial ^2v}{\partial ^2x_1}-\frac{\partial
^2v}{\partial ^2x_2}-\cdots -\frac{\partial ^2v}{\partial
^2x_m}\tag{31}
$$
Set $\Sigma =\{x\in R^m:x_1=0\}$. Obviously $\Sigma $ is spacelike
in $R^m$ with respect to the metric
$\sum_{ij}a_{ij}(0)dx_idx_j=-dx_1^2+\sum_{j=2}^mdx_j^2$.

We may write (20) briefly as follows:
$$F(\partial _i\partial
_ju)=0\tag{32}
$$
where $F=F(\zeta _{ij})$ is a polynomial in its argument $\zeta
=(\zeta _{ij})=(\partial _i\partial _ju)$, $1\leq i,j\leq m$.
Obviously
$$a^{ij}=\partial _{\zeta _{ij}}F=\frac{\partial F(\zeta
)}{\partial\zeta _{ij}}\tag{33}
$$

We prescribe the Cauchy data on $\Sigma =\{x_1=0\}$ as follows:
$$
\cases &u|_{x_1=0}=f(x^{\prime })\\
&\frac{\partial u}{\partial x_1}|_{x_1=0}=h(x^{\prime}),\quad
x^{\prime}=(x_2,...,x_m)\in R^{m-1} \endcases\tag{34}
$$
where $f,h\in C_0^\infty (R^{m-1})$ (or $S(R^{m-1})$ the Schwartz
space of rapidly decreasing functions). Then $\partial _i\partial
_ju=\partial _i\partial _jf$ on $\Sigma $ for $2\leq i,j\leq m$,
$\partial _1\partial _iu=\partial _ih$ on $\Sigma $ for $2\leq
i\leq m$. Notice that $F(\zeta _{ij})$ is a sum of terms of the
form $\zeta _{i_1i_2}\cdots \zeta _{i_kj_k}$ with $k$ odd
($k=1,3,...,[\frac m2+1]$) and $(i_s,j_s)\neq (i_t,j_t)$ if $s\neq
t$, $1\leq s,t\leq k$. So $(\partial _{\zeta _{11}}F)(\zeta
_{ij})$ is a polynomial in $\zeta _{ij}$ with $(i,j)\neq (1,1)$.
Consequently, (34) uniquely specifies $\partial _1^2u$, hence
$((\partial _{\zeta _{ij}}F)(\zeta ))$ on $\Sigma $, provided
$$
(\partial _{\zeta _{11}}F)(\partial _i\partial _ju)\neq 0\tag{35}
$$
everywhere on $\Sigma $. From Theorem 3.5 and (33), we know that
$((\partial _{\zeta _{ij}}F)(\zeta ))$ has index $1$. We may
choose the Cauchy data (34) such that $\Sigma $ is spacelike with
respect to the metric $\sum_{ij}a_{ij}dx_idx_j$, where
$(a_{ij})=(a^{ij})^{-1}$. Hence we have (cf. [Ta], [H\"or])

\proclaim{Proposition 3.6} Under the above hypothesis on the
Cauchy data (35), then the Cauchy problem (32), (35) has a unique
smooth solution on some neighborhood of $\Sigma $ in $R^m$.
\endproclaim
\remark{Remark 3.3}The Cauchy data satisfying (35) and the
spacelike condition can always be prescribed. For example, we have
already known that the trivial solution $u\equiv 0$ of (20)
corresponds to the Cauchy data $(f,h)=(0,0)$. By continuity, it is
easy to see that the Cauchy data $(f,h)$ always satisfy (35) and
the spacelike condition on $\Sigma $ if $(f,g)$ is in a small
neighborhood of $(0,0)$ in the functional space $C_0^\infty
(R^{m-1})\times C_0^\infty (R^{m-1})$.
\endremark

Finally, we consider the global existence problem for Eq. (32)
with the following Cauchy data
$$
\cases &u|_{x_1=0}=\varepsilon f\\
&\frac{\partial u}{\partial x_1}|_{x_1=0}=\varepsilon
h\endcases\tag{36} $$ where $f,h\in C_0^\infty (R^{m-1})$ and
$\varepsilon >0$ is small.

Set $f(\partial _i\partial _ju)=F(\partial _i\partial _ju)+\square
u$. Then Eq. (32) (i.e. (20)) is equivalent to
$$
\square u=f(\partial _i\partial _ju)\tag{37}
$$
where $f$ is a smooth function of $\zeta _{ij}=\partial _i\partial
_ju$ and vanishes of third order at $0$. It is known from [Kl1]
that (37) and (36) has a $C^\infty $ global solution for
sufficiently small $\varepsilon $ when $m-1\geq 4$. For $m-1=3$,
we know from [Kl2] that (37) and (36) also admits a $C^\infty $
global solution for sufficiently small $\varepsilon $, because the
Taylor expansion of $f(\zeta _{ij})$ in some neighborhood of
$(u,\partial u,\partial ^2u)=0$ does not contain any quadratic
term.

The remaining case is $m-1=2$. In this case, we get from (24) the
following
$$
f(\partial _i\partial _ju)=\det (\partial _i\partial
_ju)
$$
which is a homogeneous polynomial of degree $3$ in $\zeta
_{ij}=\partial _i\partial _ju$. Thus we meet the critical case. To
establish the global existence, we should verify the so-called
null condition, which was first introduced by Klainerman [Kl2] for
the case $m-1=3$ (see also [Chr]).

Setting $W=(W_0,W_1,W_2,W_3)=(u,\partial _1u,\partial _2u,\partial
_3u)$, we find that solving (37) is equivalent to solve
$$
\cases &\square W_0=det(\partial _iW_j)\\
&\square W_1=C(\partial W)^{ij}\partial _i\partial _jW_1\\
&\square W_2=C(\partial W)^{ij}\partial _i\partial
_jW_2\\
&\square W_3=C(\partial W)^{ij}\partial _i\partial
_jW_3\endcases\tag{38}
$$
where $(C(\partial W)^{ij})$ is the cofactor matrix of $(\partial
_iW_j)$ $(1\leq i,j,k\leq 3)$.

\definition {Definition 3.1}Let $G=G((w_a);(\partial _iw_b);(\partial
_{ij}^2w_c))$ be a smooth function of $w_a$ ($a=0,1,...,N$),
$\partial _iw_a$($a=0,1,...,N$) and $\partial
_{ij}^2w_b$($a=0,1,...,N$, $i,j=1,...,m$). We say that $G$
satisfies the null condition when
$$
G((\lambda _a);(\mu _bX_i);(v_cX_jX_k))=0
$$
for all $\lambda ,\mu ,v\in R^{N+1}$ and all $X=(X_1,\cdots
,X_m)\in R^m$ satisfying $X_1^2-X_2^2-\cdots -X_m^2=0$.
\enddefinition

We may verify directly that the functions $\det (\partial _iW_j)$,
$C(\partial W)^{ij}\partial _i\partial _jW_k$ ($k=1,2,3$)
appearing on the right hand side of (38) satisfy the null
condition in Definition 3.1. Consequently, we know from Theorem
1.2 of [Ka] that the Cauchy problem (37),(36) has a unique global
$C^\infty -$solution too when $m-1=2$.

In conclusion, we have shown the following

\proclaim{Theorem 3.7}The Cauchy problem (37), (36) with $f,h\in
C_0^\infty (R^{m-1})$ has a unique $C^\infty $ global solution $u$
if $m\geq 3$ and $\varepsilon $ is sufficiently small.
\endproclaim
\remark{Remark 3.4}For a solution $u$ of (20), the graph of
$(\nabla u)I_{k,m}$ is an indefinite special Lagrangian
submanifold provided that (19) is satisfied. When $\varepsilon $
is sufficiently small, the solution obviously satisfies the
non-degenerate condition (19) everywhere.
\endremark

Proposition 3.3, Proposition 3.6 and Theorem 3.7 show that
timelike special Lagrangian submanifolds exist in abundance. We
will construct more nontrivial explicit examples of indefinite
special Lagrangian submanifolds in next section.

\heading{\bf 4. Explicit examples of indefinite special Lagrangian
submanifolds}
\endheading
\vskip 0.3 true cm

In this section, we hope to construct some explicit indefinite
special Lagrangian submanifolds by the following two methods: the
moment map method for symmetric indefinite special Lagrangian
submanifolds and the normal bundle constructions.

\vskip 0.3true cm {\bf 4.1 Symmetric indefinite special Lagrangian
submanifolds} \vskip 0.3 true cm Let $G$ be a connected Lie group
of holomorphic isometries of $C_k^m$. Let $g $ be the Lie algebra
of $G$, and $g^{*}$ the dual space of $g$. Then a moment map for
the $G$-action on $C_k^m$ is a smooth map $\mu $ $:C^m\rightarrow
g^{*}$ such that (a) $d(\mu ,\xi )=i_{X_\xi }\omega $ for all $\xi
\in g$, where $(\ ,\ )$ denotes the pairing between $g^{*}$ and
$g$, and $X_\xi $ is the infinitesimal action corresponding to
$\xi $; (b) $\mu (kx)=$ $Ad_{k^{-1}}^{*}\mu (x),$ $\forall k\in G$
and $x$ $\in M$, where $Ad$ denotes the coadjoint action (For
basic properties of the moment maps, the reader could refer to
[Si]).

According to the terminology of Symplectic geometry, a $G-$action
is called Hamiltonian if it admits a moment map. By using the
properties (a) and (b) of a moment map, it is easy to prove the
following:

\proclaim{Lemma 4.1}(Cf. [Jo2]) Let $G\times M\rightarrow M$ be a
Hamiltonian action on a symplectic manifold $M$. If $L$ is a
connected $G$-invariant Lagrangian submanifold in $M$, then
$M\subset \mu ^{-1}(c)$ for some $c$ $\in $ $Z(g^{*})$, where
$Z(g^{*})$ denotes the center of $g^{*}$.
\endproclaim

First, let's determine the moment map of the natural action of
$SU(k,m-k)$ on $C_k^m$ : $(A,z)\longmapsto Az$. Its infinitesimal
action is given by
$$
X_\xi (z)=\xi z
$$
where $\xi \in u(k,m-k)$ and $z\in C^m$ is a column vector. We fix
the following inner product on $u(k,m-k)$,
$$
<\xi ,\eta >:=-tr(\xi \eta)\tag{39}
$$
to identify $u(k,m-k)$ with $u^{*}(k,m-k)$. Two vectors $v,w\in
C_k^m$, induce a complex linear map
$$
vw^{*}I_{k,m}:C^m\rightarrow C^m
$$
where $w^{*}=(\overline{w}_1,...,\overline{w}_m)$ is the conjugate
transpose of $w$. Obviously we have
$$
vw^{*}I_{k,m}(z)=h_{k,m}(z,w)v\tag{40}
$$
The symplectic structure associated to the inner product $h_{k,m}$
is $$\omega _{(k,m)}=-Im(h_{k,m})$$

\proclaim{Proposition 4.2}The action of $SU(k,m)$ on $C_k^m$ is
Hamiltonian with moment map
$$
\mu (z)=-\frac i2zz^{*}I_{k,m}\tag{41}
$$
\endproclaim
\demo{Proof}For $\xi \in u(k,m-k)$ and $v,w\in C_k^m$, we have
$$
\aligned tr[\xi \cdot \frac i2(vw^{*}I_{k,m}+wv^{*}I_{k,m})
&=\frac i2h_{k,m}(\xi v,w)+\frac i2h_{k,m}(\xi w,v)\\
&=\frac i2h_{k,m}(\xi v,w)-\frac i2h_{k,m}(w,\xi v)\\
&=\frac i2h_{k,m}(\xi v,w)-\frac i2\overline{h_{k,m}(\xi v,w)}\\
&=-Im(h_{k,m}(\xi v,w))\\
&=\omega _{(k,m)}(\xi v,w)\endaligned
$$
The defining equation for a moment map is
$$
d<\mu (z),\xi >(v)=\omega _{(k,m)}(X_\xi (z),v)=\omega
_{(k,m)}(\xi z,v)
$$
i.e.,
$$
<d\mu (z)v,\xi >=tr[\xi \cdot \frac
i2(zv^{*}I_{k,m}+vz^{*}I_{k,m})]\tag{42}
$$
Using the inner-product on $u(k,m-k)$, we get from (42) the
following
$$
-tr[d\mu (z)v\xi ]=tr[\xi \cdot \frac
i2(zv^{*}I_{k,m}+vz^{*}I_{k,m})]
$$
which is satisfied by
$$
\mu (z)=-\frac i2zz^{*}I_{k,m}\tag{43}
$$
It is easy to verify that the map $\mu $ given by (43) satisfies
the equivariant property. Therefore we prove this proposition.\qed

\enddemo

First we hope to construct some $T^{m-1}-$invariant indefinite
special Lagrangian submanifolds. Here $T^{m-1}$ is the subgroup
$$
T^{m-1}=\{diag(e^{i\theta _1},...,e^{i\theta _m}):\theta _1+\cdots
+\theta _m=0\}\tag{44}
$$
in $SU(k,m-k)$.

\proclaim{Lemma 4.3}Suppose $T^{m-1}$acts on $C_k^m$ as follows:
$$
(e^{i\theta _1},...,e^{i\theta _{m-1}})\cdot z=\left(\matrix
e^{i\theta _1}z_1\\
\vdots\\
e^{i\theta_{m-1}}z_{m-1}\\
e^{i\theta _m}z_m\endmatrix\right )
$$
where $\theta _m=-\theta _1-\cdots -\theta _{m-1}$. Then the
moment map of this action is given by (up to a
constant)
$$
\mu(z)=diag(|z_1|^2+|z_m|^2,...,|z_k|^2+|z_m|^2,
|z_m|^2-|z_{k+1}|^2,...,|z_m|^2-|z_{m-1}|^2)
$$
\endproclaim
\demo{Proof}We consider the homomorphism $\varphi
:T^{m-1}\rightarrow U(k,m-k)$ defined by
$$
\varphi(diag(e^{i\theta _1},...,e^{i\theta
_{m-1}}))=diag(e^{i\theta _1},...,e^{i\theta
_m})\tag{45}
$$
where $\theta _m=-\theta _1-\cdots -\theta _{m-1}$. Then the
induced homomorphism $d\varphi :t\cong R^{m-1}\rightarrow
su(k,m-k)$ between their Lie algebras gives
$$
\aligned d\varphi (\frac\partial {\partial \theta _1})&=i\
diag(1,....,-1)\\
&\cdots \\
d\varphi (\frac \partial {\partial \theta _{m-1}})&=i\
diag(0,...,1,-1)\endaligned\tag{46}
$$
Therefore, using the inner product on $u(k,m-k)$, we have
$$
<\mu ,d\varphi (\frac \partial {\partial \theta
_j})>=\cases \frac i2(|z_j|^2+|z_m|^2),\quad j=1,...,k\\
-\frac i2(|z_j|^2-|z_m|^2),\quad j=k+1,...,m-1\endcases\tag{47}
$$
Using the inner product, we may regard the moment map of the
$T^{m-1}-$action as a map into its Lie algebra. Therefore we prove
the Lemma.\qed
\enddemo

\proclaim{Theorem 4.4}Let $F=(f_1,...,f_{m-1},f_m):C^m\rightarrow
R^m$ be a map defined by
$$
f_j=\cases |z_j|^2+|z_m|^2,\quad j=1,...,k\\
|z_j|^2-|z_m|^2,\quad j=k+1,...,m-1\endcases
$$
and
$$
f_m=\cases Re(z_1...z_m),\quad \text{if $m$ is
even,}\\
Im(z_1...z_m),\quad \text{if $m$ is odd}\endcases
$$
where $m\geq 3$. Set $D=\{z\in C^m:\det (\frac{\partial
f_i}{\partial \overline{z}_j})=0\}$. Let $M_c=F^{-1}(c)$ be the
inverse image of a point $c\in (R^{+})^k\times R^{m-k}-F(D)$,
where $R^{+}=\{x\in R:x>0\}$. Then $M_c$ (with the correct
orientation) is an indefinite special Lagrangian submanifolds of
$C_k^m$.
\endproclaim
\demo{Proof}First we assume that $m$ is even, i.e.,
$$
f_m(z)=\frac{z_1...z_m+\overline{z}_1...\overline{z}_m}2\tag{48}
$$
$$
f_i=\cases z_i\overline{z}_i+z_m\overline{z}_m,\quad i=1,...,k\\
z_i\overline{z}_i-z_m\overline{z}_m,\quad i=k+1,...,m-1
\endcases\tag{49}
$$
A direct computation shows that
$$
\det(\frac{\partial f_i}{\partial \overline{z}_j})=\frac
12\{-\sum_{i=1}^k|z_1|^2\cdots \widehat{|z_i|^2}\cdots
|z_m|^2+\sum_{j=k+1}^m|z_1|^2\cdots \widehat{|z_j|^2}\cdots
|z_m|^2\}\tag{50}
$$
and
$$
\overline{z}_i\frac{\partial f_m}{\partial
\overline{z}_i}=\overline{z}_m\frac{\partial f_m}{\partial
\overline{z}_m},\quad 1\leq i\leq m-1\tag{51}
$$
Note that $\det (\frac{\partial f_i}{\partial \overline{z}_j})$ is
a real analytic function on $C^m$. So the critical points set of
$F=(f_1,...,f_m):C^m\rightarrow R^m$ is a real hypersurface of
$C^m$ given by
$$
D=\{z\in C^m:\det (\frac{\partial f_i}{\partial
\overline{z}_j})=0\}\tag{52}
$$
Its image $F(D)$ has measure zero for the usual measure on $R^m$.
Thus, for any $c=(c_1,...,c_m)\in (R^{+})^k\times R^{m-k}-F(D)$,
$M_c$ is a smooth manifold of dimensional $m$. The condition (16)
follows from (51) and the property of the moment map. From the
proof of Proposition 2.4, we see that $\nabla ^gf_i$ corresponds
to the complex vector
$$
\aligned &2(-\frac{\partial f_i}{\partial \overline{z}_1},\cdots
,-\frac{\partial f_i}{\partial \overline{z}_k},\frac{\partial
f_i}{\partial \overline{z}_{k+1}},\cdots ,\frac{\partial
f_i}{\partial \overline{z}_m})\\
=&\cases
2(0,...,0,-z_i,0,...,z_m)\quad\text{if $1\leq i\leq k$,}\\
2(0,...,0,z_i,0,...,-z_m)\quad\text{if $k+1\leq i\leq m-1$,}\\
(-\overline{z}_2\cdots \overline{z}_m,...,-\overline{z}_1\cdots
\widehat{\overline{z}}_k\cdots \overline{z}_m,\overline{z}_1\cdots
\widehat{\overline{z}}_{k+1}\cdots
\overline{z}_m,...,\overline{z}_1\cdots \widehat{\overline{z}}_m),
&\text{if $i=m$}\endcases\endaligned
$$
under the identification $R^{2m}$ with $C^m$. It follows
that
$$
\sum_{l=1}^m\varepsilon _l\frac{\partial f_i}{\partial
\overline{z}_l}\frac{\partial f_j}{\partial z_l}\in R,\quad
\text{for $1\leq i,j\leq m$}
$$
Since $\{f_i\}$ are real valued, we get
$$
(g_{(2k,2m)}(\nabla ^gf_i,\nabla ^gf_j))=4(\partial f_i/\partial
\overline{z}_l)I_{k,m}(\partial f_i/\partial z_l)^t\tag{53}
$$
This implies that $\det (g_{(2k,2m)}(\nabla ^gf_i,\nabla
^gf_j))\neq 0$ if and only if $\det (\partial f_i/\partial
\overline{z}_l)\neq 0$. Hence for any $c\in (R^{+})^k\times
R^{m-k}-F(D)$, $\{f_i\}$ satisfies the conditions of Lemma 2.3 on
$M_c$. Therefore $M_c$ is a Lagrangian submanifold in $C^m$ with
respect to $\omega _{(k,m)}$. From the expression of $f_m$, we get
$$
\overline{z}_m\frac{\partial f_m}{\partial
\overline{z}_m}z_1\cdots z_m=\frac 12|z_1\cdots z_m|^2\tag{54}
$$
for $i=1,...,m$. Obviously $\det_C(\partial f_i/\partial
\overline{z}_j)$ is a sum of terms of the form
$$
\pm \frac{\overline{z}_l}{|z_l|^2}\frac{\partial f_m}{\partial
\overline{z}_l}z_1\cdots z_m=\pm
\frac{\overline{z}_m}{|z_l|^2}\frac{\partial f_m}{\partial
\overline{z}_m}z_1\cdots z_m\tag{55}
$$
Consequently
$$
Im\det (\partial f_i/\partial \overline{z}_j)=0
$$
By Proposition 2.4, we get the result for the case of $m$ even. We
may prove the similar result for the case of $m$ odd. \qed
\enddemo

\remark{Remark 3.1}If $c\in F(D)$, $M_c$ may not be a Lagrangian
submanifold or have various kinds of singularity depending on $c$.
\endremark

Next, we consider the subgroup $i:SO(k,m-k)\rightarrow SU(k,m-k)$
which acts diagonally on $C^m\cong R^m\oplus R^m$. The Lie algebra
of $SO(k,m-k)$ is
$$
so(k,m-k)=\{A\in M(m,R):A^tI_{k,m}+I_{k,m}A=0\}\tag{56}
$$
Using the natural basis of $so(k,m-k)$ and the moment map of
$SU(k,m-k)$, we may obtain the moment map of $SO(k,m-k)-$action as
follows:
$$
\aligned \mu (z)=&\sum_{1\leq i<j\leq
k}Im(z_i\overline{z}_j)(E_{ij}-E_{ji})+\sum_{{1\leq i\leq k}, {k+1\leq j\leq m}}Im(z_i\overline{z}_j)(E_{ij}+E_{ji})\\
&-\sum_{k+1\leq i<j\leq
m}Im(z_i\overline{z}_j)(E_{ij}-E_{ji})\endaligned\tag{57}
$$
where $E_{ij}$ denotes the $m\times m$ matrix such that the entry
at the $i-$th row and $j-$th column is $1$ and other entries are
all zero.

As $Z(g^{*})=0$, any $SO(k,m-k)-$invariant indefinite Lagrangian
$m-$fold lies in $\mu ^{-1}(0)$. Now every point in $\mu ^{-1}(0)$
may be written as $(\lambda t_1,...,\lambda t_m)$, where $\lambda
\in C$ and $t:=(t_1,...,t_m)\in R_k^m$ is normalized so that
$$
\sum \varepsilon _jt_j^2=\cases 1 &\text{if $t$ is
spacelike},\\
0 &\text{if $t$ is lightlike},\\
-1 &\text{if $t$ is timelike}\endcases\tag{58}
$$
If $(z_1,...,z_m)\in \mu ^{-1}(0)$ satisfies
$-\sum_{j=1}^kz_j\overline{z}_j+\sum_{j=k+1}^mz_j\overline{z}_j>0$
(resp. $<0$), then the $SO(k,m-k)-$orbit $\Theta _z\cong
S_k^{m-1}(r^2)$ the pseudo-Riemannian sphere (resp.
$H_{k-1}^{m-1}(-r^2)$ the pseudo-hyperbolic space). First we note
that, for any regular curve $\Gamma\subset C^{*}=C\backslash
\{0\}$, the submanifold defined by
$$
M_\Gamma =\{(x,y)=\lambda (t_1,...,t_m)\in
C_k^m:\sum_{j=1}^m\varepsilon _jt_j^2=\pm 1,\ \lambda \in \Gamma
\}\tag{59}
$$
is Lagrangian w.r.t. $\omega _{(k,m)}$. In fact, we may write
$\lambda=\xi (s)+i\eta (s)$ and compute the induced $2-$form of
$\omega _{(k,m)}$ on $M$ as follows:
$$
\aligned \omega _{(k,m)}|_{M_\Gamma } &=\sum_{j=1}^m\varepsilon
_jdx_j\wedge dy_j\\&=\sum_{j=1}^m\varepsilon _jd(\xi t_j)\wedge
d(\eta t_j)\\
&=\sum_{j=1}^m\varepsilon _j(t_jd\xi +\xi dt_j)\wedge (t_jd\eta
+\eta dt_j)\\&=\sum_{j=1}^m\varepsilon _jt_j^2d\xi \wedge d\eta
+\sum_{j=1}^m\varepsilon _j(t_j\eta d\xi \wedge dt_j+t_j\xi
dt_j\wedge d\eta)\\
&=\pm d\xi \wedge d\eta\endaligned
$$
where we use the fact $\sum_j\varepsilon _jt_j^2=\pm 1$ in the
last equality. The fact $\dim _R\Gamma =1$ leads to $\omega
_{(k,m)}|_{M_\Gamma }\equiv 0$. The induced metric on $M_\Gamma $
is given by
$$
ds_{M_\Gamma }^2=\pm |\lambda ^{\prime }|^2ds^2+\lambda
^2dt^2
$$
where $dt^2$ is the metric of $S_k^{m-1}(1)$ or
$H_{k-1}^{m-1}(-1)$. It is easy to see that $ds_{M_\Gamma }^2$ is
non-degenerate.

\proclaim{Theorem 4.5}Let
$$
M_c=\{(x,y)=\lambda (t_1,...,t_m)\in C_k^m:\sum_{j=1}^m\varepsilon
_jt_j^2=1,\ \lambda \in C,\ Im(\lambda ^m)=c\}
$$
or
$$
M_c=\{(x,y)=\lambda (t_1,...,t_m)\in C_k^m:\sum_{j=1}^m\varepsilon
_jt_j^2=-1,\ \lambda \in C,Im(\lambda ^m)=c\}
$$
Then $M_c$ (with the correct orientation) is an indefinite special
Lagrangian submanifold of $C_k^m$.
\endproclaim
\demo{Proof}From the above discussion, we have already known that
$M_c$ is a Lagrangian submanifold. Locally we may express $\eta $
as a function of $\xi $, say, $\eta =\varphi (\xi )$. So $M_c$ is
the graph of
$$
f(x)=\varphi (\xi )\frac x\xi\tag{60}
$$
where $\xi =\sqrt{\pm \sum_i\varepsilon _ix_i^2}$ . The
differential $\varphi _{*}$ of this map from $R^m$ to $R^m$ is
given by the matrix $(h_{ij})$ where
$$
\aligned h_{ij}&=\frac \partial {\partial x_i}[\varphi (\xi
)\frac{x_j}\xi ]\\
&=\frac{\varphi (\xi )}\xi \delta _{ij}+\frac d{d\xi
}(\frac{\varphi (\xi )}\xi )\{\pm \varepsilon _i\frac{x_ix_j}\xi
\}\endaligned\tag{61}
$$
Then the linear map $\varphi _{*}:R^m\rightarrow R^m$ has the
eigenvector $x$ with eigenvalue
$$\frac{\varphi (\xi )}\xi +\xi
\frac d{d\xi }(\frac{\varphi (\xi )}\xi )=\frac{d\varphi (\xi
)}{d\xi }\tag{62}
$$
Moreover, the hyperplane perpendicular to $x$ is an eigenspace
with eigenvalue $\frac{\varphi (\xi )}\xi $ of multiplicity $m-1$.
Set $K=I+i(h_{ij})$. Hence the graph of $f$ is special Lagrangian
if and only if
$$
Im\{\det K\}=0
$$
i.e.,
$$
\frac 1{\xi ^{m-1}}Im\{(\xi +i\eta )^{m-1}(d\xi +id\eta
)=0\tag{63}
$$
Therefore the integral curves of the O.D.E. are of the form
$Im(\xi +i\eta )^m=c$ for some $c\in R$.\qed
\enddemo
\remark{Remark 4.2}(1) For $c\neq 0$, each component of the
manifold $M_c$ is diffeomorphic to $R\times S_k^{m-1}(1)$ or
$R\times $ $H_{k-1}^{m-1}(-1)$. When $c=0$, it is a singular union
of $m$ copies of $m-$dimensional Lagrangian cones, and the link of
each cone is $S_k^{m-1}(1)$ or $H_{k-1}^{m-1}(-1).$
\newline(2) In [Ch], B. Chen introduced the notion of complex
extensors in $C_k^m$ to construct $SO(k,m-k)-$invariant Lagrangian
submanifolds. He also got the representations (59) for these
Lagrangian submanifolds, and then gave the classification of
Lagrangian $H-$umbilical submanifolds.
\endremark
\vskip 0.3 true cm {\bf 4.2 Indefinite special Lagrangian normal
bundles}\vskip 0.3 true cm

Let $M^n$ be an indefinite submanifold in the pseudo-Euclidean
space $R_k^m$. Let $g$ be the induced pseudo-Riemannian metric on
$M$ with index $s$. For a normal vector field $\xi \in \Gamma
(T^{\bot }M)$, the formula
$$
\det(tI_m-A_\xi )=\sum_{l=0}^m(-1)^l\sigma _l(\xi
)t^{m-l}
$$
defines a sequence $\sigma _l(\xi )$ of smooth functions on $M$.
Clearly, $\sigma _0(\xi )=1$ while $\sigma _1(\xi )=tr(A_\xi )$.

\definition{Definition 4.1}An indefinite submanifold $M^n$ of $R_k^m$ is said to
be austere if, for every $\xi \in \Gamma (T^{\bot }M)$ and every
integer $l$ satisfying $0\leq l\leq m/2$, we have $\sigma
_{2l+1}(\xi )=0$.
\enddefinition

Notice that a selfadjoint matrix $A$ with respect to a
pseudo-Euclidean metric may have no real eigenvalues(see page 273
of [Gr] for such an example). If $A_\xi $ does have $m$ real
eigenvalues $\lambda _1,...,\lambda _m$ for each normal vector
$\xi $, then the austere condition is equivalent to the condition
that the set of eigenvalues of $A_\xi $ is of the form
$$
(\lambda _1,...,\lambda
_m)=(a,-a,b,-b,...,c,-c,0,...,0)
$$
When $M$ is spacelike, i.e., $s=0$, it is known that $A_\xi $
always has $m$ real eigenvalues for each normal vector $\xi $, and
thus it is diagonalizable. For general case, we have the following
criteria for diagonalization:

\proclaim{Lemma 4.6}Let $M^n$ be an indefinite submanifold of
dimension $n\geq 3$ in $R_k^m$. Suppose $M$ satisfies
$$
g(A_\xi X,X)+g(X,X)>0\tag{64}
$$
for each nonzero normal vector $\xi $ and every nonzero tangent
vector $X$ at the same point. Then there exists a Lorentz basis
$\{e_v\}_{v=1}^n$ with respect to $g$ at each point such that
$$
A_\xi e_i=\lambda _ve_i\quad i=1,...,n
$$
\endproclaim
\demo{Proof} The result follows immediately from Theorem 9.11 in
[Gr].\qed\enddemo \remark{Remark 4.3}A spacelike submanifold
automatically satisfies the condition (64) and the assumption
$n\geq 3$ is not necessary for the spacelike case. However, the
assumption $n\geq 3$ is necessary for the general case.\endremark

Now we assume that $M$ is a spacelike submanifold or an indefinite
submanifold satisfying Lemma 4.6. We define the embedding
$$
\psi :T^{\bot }M\rightarrow R_k^m\oplus R_k^m=C_k^m\tag{65}
$$
by setting $\psi (v_x)=(x,v(x))$ where the second factor $v(x)$ is
a vector based at the origin obtained by moving $v_x$ to the
origin. Near $x_0$ we choose a Lorentz tangent frame field
$e_1,...,e_n$ and a Lorentz normal frame field $v_1,...,v_p$,
$n+p=m$, such that $(e_1,...,v_p)$ is positively oriented and
$(\nabla ^{\perp }v_i)_{x_0}=0$.

Obviously the tangent space to this embedding at
$v(x_0)=\sum_jc_jv_j$ is spanned by the vectors
$$
\aligned
E_j&=\psi _{*}(e_j)=(e_j,A_ve_j),\quad j=1,...,n\\
N_j&=\psi _{*}(\partial /\partial t_j)=(0,v_j),\quad
j=n+1,...,m\endaligned\tag{66}
$$
It is easy to see from (9) and (66) that
$$
g_{(2k,2m)}(JN_j,N_k)=g_{(2k,2m)}(JN_j,E_l)=-g_{(2k,2m)}(N_j,JE_l)=0
$$
for all $j,k,l$. Moreover
$$
\aligned g_{(2k,2m)}(JE_j,E_k)
&=g_{(2k,2m)}((-A_ve_j,e_j),(e_k,A_ve_k)\\
&=-g_{(2k,2m)}(A_ve_j,e_k)+g_{(2k,2m)}(e_j,A_ve_k)\\
&=0
\endaligned
$$
Obviously the induced metric on $\psi (T^{\bot }M)$ from
$g_{(2k,2m)}$ is non-degenerate. Hence $\psi (T^{\bot }M)$ is a
Lagrangian submanifold of $C_k^m$ with respect to $\omega
_{(k,m)}$.

The hypothesis about $M$ implies that we may choose a Lorentz
basis $e_1,...,e_n\ $at $x_0$ such that $A_v(e_j)=\lambda _je_j$,
$j=1,...,n$. Consequently, up to a sign, the tangent plane
$\varsigma $ of the embedding $\psi $ at $v_{x_0}$ is given by
$$
\aligned \varsigma &=E_1\wedge \cdots \wedge E_n\wedge N_1\wedge
\cdots \wedge E_p\\
&=(e_1,\lambda _1e_1)\wedge \cdots \wedge (e_n,\lambda
_ne_n)\wedge (0,v_1)\wedge \cdots \wedge (0,v_p)\endaligned
\tag{67}
$$
Since $e_1,....,e_n,v_1,...,v_p$ is a Lorentz basis, we may
reorder the basis and then perform an $SO(k,m-k)$ change of
coordinates on $R_k^m$ such that $\{e_1,...,e_n,v_1,...,v_p\}$
becomes the standard Lorentz basis of $R_k^m$. It follows that
$$
dz_1\wedge \cdots \wedge dz_m(\varsigma
)=i^p\prod_{j=1}^n(1+i\lambda _j)\tag{68}
$$

\proclaim{Theorem 4.7}Let $M^n$ be a space-like submanifold or an
indefinite submanifold satisfying Lemma 3.7 in $R_k^m$. Then the
normal bundle $\psi (T^{\bot }M)$ is indefinite special Lagrangian
in $C_k^m=R_k^m\oplus R_k^m$ if and only if $M$ is austere in
$R_k^m$.\endproclaim

The above result shows that it would be interesting to find more
austere submanifolds. Maximal spacelike surfaces in $R_1^3$ are
automatically austere. The Weierstrass formula in [Ko] provides us
many examples of maximal surfaces. By generalizing Bryant's idea
in [Br], the authors in [DH] also construct some examples of
spacelike austere submanifolds in pseudo-Euclidean spaces of
higher dimensions.
\vskip 0.3 true cm

{\bf Acknowledgments}: The author would like to thank Professor B.
Chen for sending him related papers. He would also like to thank
Professors J. Hong, N. Mok, Y. Zhou and Dr. Y. Han for their
helpful discussions.

\heading {\bf Appendix: Instability of indefinite minimal
submanifolds}\endheading

In this appendix, we will investigate the stability problem of an
$m-$dimensional indefinite minimal submanifold $\varphi
:M\rightarrow R_{\ n}^N$. For simplicity, we sometimes write the
metric of $R_n^N$ as $<\ ,\ >$. Suppose that the induced
pseudo-Riemannian metric $g$ on $M$ has index $k$ with $0<k<m$.
Using a coordinate system $(u^i)$ of $M$, $g$ is expressed as
$$
g=g_{ij}du^idu^j\tag{69}
$$
The volume functional for $\varphi $ is defined by
$$
V(\varphi )=\int_M\sqrt{(-1)^k\det (g_{ij})}du^1\wedge \cdots
\wedge du^m\tag{70}
$$
For a variation $\varphi _t$ corresponding to a normal vector
field $W$ along $M$ with compact support, we set $V(t)=V(\varphi
_t)$. It is known that the submanifold $M$ is minimal for the
functional $V$ if and only if $H\equiv 0$.

Now we assume that $M$ is an indefinite minimal submanifolds of
$R_n^N$. Set $dv=\sqrt{(-1)^k\det (g_{ij})}du^1\wedge \cdots
\wedge du^m$. By the usual computation, we have the second
variational formula:
$$
V^{\prime \prime }(0)=\int_M\{<\nabla ^{\bot }W,\nabla ^{\bot
}W>-<h\circ h^t(W),W>\}dv\tag{71}
$$
where $h$ is the second fundamental form and
$$
<h\circ h^t(W),W>=g^{ik}g^{jl}<h(\frac
\partial {\partial u^i},\frac
\partial {\partial u^j}),W><h(\frac \partial {\partial u^k},\frac
\partial {\partial u^l}),W>\tag{72}
$$

Choose a local spacelike vector field $X$ with $|X|=1$ on $\Omega
\subset M$. By the Flow-box theorem, there exists a cubic
coordinate system $(D,\psi ,u^i)$ $\subset \subset \Omega $ with
$$
\psi (D)=\{(u^1,...,u^m):-\delta <u^i<\delta \}\tag{73}
$$
such that
$$
X|_D=\frac \partial {\partial u^1}\tag{74}
$$
We may choose $\delta $ sufficiently small so that the normal
bundle is trivialized on $D$. Thus there exists a local normal
vector field $\xi $ on $D$ with the property that $<\xi ,\xi >\neq
0$ everywhere on $D$. Without lose of generality, we may assume
that $\xi $ is spacelike everywhere on $D$. Set
$$
W=f(u^1,...,u^m)\xi\tag{75}
$$
where
$$
f(u^1,...,u^m)=[1+\cos \frac{(2q+1)\pi }\delta u^1]\rho
(u^2,...,u^m)
$$
where $q$ is an integer and $\rho \in C_c^\infty (\Omega )$. So
$W|_{\partial \Omega }=0$. Using $W$ as a variation vector field,
we get from (71) the following:
$$
\aligned V^{\prime \prime }(0)&=\int_M\{<\nabla _{\frac \partial
{\partial u^1}}^{\bot }f\xi ,\nabla _{\frac \partial {\partial
u^1}}^{\bot }f\xi >+2\sum_{j=2}^mg^{1j}<\nabla _{\frac \partial
{\partial u^1}}^{\bot }f\xi ,\nabla _{\frac \partial {\partial
u^j}}^{\bot }f\xi >\\
&+\sum_{j,l=2}^mg^{jl}<\nabla _{\frac \partial {\partial
u^j}}^{\bot }f\xi ,\nabla _{\frac \partial {\partial u^l}}^{\bot
}f\xi >-f^2<h\circ h^t(\xi ),\xi >\}dv\\
&=\frac{(2q+1)^2\pi ^2}{\delta ^2}\int_D\rho ^2\sin
^2(\frac{(2q+1)\pi }\delta u^1)dv+\frac{(2q+1)}\delta
I_1+I_2\endaligned\tag{76}
$$
where
$$
|I_1|\leq C_1,\quad |I_2|\leq C_2 \tag{77}
$$
Here $C_1$ and $C_2$ are constants independent of $q$. So
$$
V^{\prime \prime }(0)\geq \frac{(2q+1)^2\pi ^2}{\delta
^2}\int_D\rho ^2\sin ^2(\frac{(2q+1)\pi }\delta
u^1)dv-C_1\frac{(2q+1)}\delta -C_2\tag{78}
$$
Notice that $\int_D\rho ^2\sin ^2(\frac{(2q+1)\pi }\delta u^1)dv$
increases as $q\rightarrow +\infty $. Thus
$$
\int_D\rho ^2\sin ^2(\frac{(2q+1)\pi }\delta u^1)dv\geq
C_0>0\tag{81}
$$
where $C_0$ is a constant independent of $q$. By choosing a
sufficiently large $q$, we have
$$
V^{\prime \prime }(0)>0\tag{82} $$ As a result, the variation
increases the volume of $M$.

Similarly we may start with a timelike vector field $Y$ with
$<Y,Y>=-1$ and choose a cubic coordinate system $(u^i)$ with
$Y=\frac \partial {\partial u^1}$. Using the same variation vector
field $W$, we may get $V^{\prime \prime }(0)<0$ for a sufficiently
large $q$. In conclusion, we have proved the following:

\proclaim{Theorem A}Let $M$ be an indefinite minimal submanifold
of $R_n^N$ with index $0<k<m$. Then for any domain on $M$, there
exists a smooth variation with fixed boundary that increases the
volume, and there exists a smooth variation that decreases the
volume.\endproclaim

\remark{Remark}Such kind of instability of indefinite minimal
submanifolds was first obtained by Gorokh [Go] for timelike
minimal surfaces in $R_1^3$. Here we generalize his result to the
case of any dimension and codimension. As a consequence, we know
that there is no minimizing property or maximizing property for
the indefinite special Lagrangian submanifolds.\endremark

\vskip 1 true cm

\vskip 1 true cm

Institute of Mathematics

Fudan University, Shanghai 200433

P.R. China

And

Key Laboratory of Mathematics

for Nonlinear Sciences

Ministry of Education

\vskip 0.2 true cm
yxdong\@fudan.edu.cn


\Refs \widestnumber\key{SW1}

\ref\key AAI \by P. Allen, L. Andersson, J. Isenberg\paper
Timelike minimal submanifolds of general co-dimension in Minkowski
spacetime\paperinfo  math.DG/0512036
\endref

\ref\key Bre\by S. Brendle\paper Hypersurfaces in Minkowski space
with vanishing mean curvature\paperinfo Comm. Pure Appl. Math., 55
(2002), 1249-1279\endref

\ref\key Br\by R. Bryant\paper Some remarks on the geometry of
austere manifolds\paperinfo Bol. Soc. Brasil. Mat. 21, 122-157
(1991)\endref

\ref\key Ch\by B.Y. Chen\paper Complex extensors and Lagrangian
submanifolds in indefinite complex Euclidean spaces\paperinfo
Bull. of the Inst. of Math. Acad. Sinica Vol. 31, No. 3,
2003\endref

\ref\key Chr\by D. Christodoulou\paper Global solutions of
nonlinear hyperbolic equations for small initial data\paperinfo
Comm. Pure Appl. Math. 39 (1986), 267-282\endref

\ref\key De1\by T. Deck\paper A geometric Cauchy problem for
timelike minimal surfaces \paperinfo Ann. Global Anal. Geom. 12
(1994), 305-312\endref

\ref\key De2\by T. Deck\paper Classical string dynamics with
non-trivial topology\paperinfo J. Geom. Phys. 16 (1995),
1-14\endref

\ref\key DH\by Y.X. Dong, Y.B. Han\paper On spacelike Austere
submanifolds\paperinfo to appear\endref

\ref\key Go\by V.P. Gorokh\paper The stability of a minimal
surface in a pseudo-Euclidean space\paperinfo Journal of
Mathematical Sciences, Volume 53, Number 5, (1991), 491-493\endref

\ref\key Gr1\by W. Greub\paper Linear Algebra\paperinfo Graduate
Text in Math. 23, Springer-Verlag, 1975\endref

\ref\key Gu1\by C.H. Gu\paper On the motion of a string in a
curved space-time\paperinfo in: Hu Ning (ed), Proc. 3rd. Grossmann
meeting on General Relativity, Beijing, 139-142, (1983)\endref

\ref\key Gu2\by C.H. Gu\paper A global study of extremal surfaces
in 3-dimensional Minkowski space\paperinfo in: ChaoHao Gu, Marcel
Berger and Robert L. Bryant (eds.), Differential Geometry and
Differential Equations (Shanghai 1985), Lecture Notes in
Mathematics, Vol. 1255, Springer-Berlin, (1985), 173-180\endref

\ref\key Gu3\by C. H. Gu\paper The extremal surfaces in the
3-dimensional Minkowski space\paperinfo Acta Math. Sinica, New
Series 1 (1985), 173-180\endref

\ref\key Hi\by N.J. Hitchin\paper The moduli space of special
Lagrangian submanifolds\paperinfo Dedicated to Ennio De Giorgi,
Ann. Scuola Norm. Sup. Pisa Cl. Sci. (4) 25 (1997), no. 3-4,
503-515 (1998)\endref

\ref\key HL\by R. Harvey, H.B. Lasson. Jr.\paper Calibraed
Geometries\paperinfo Acta Math. 148 (1982), 47-157\endref

\ref\key H\"o\by L. H\"ormander\paper Lectures on nonlinear
hyperbolic differential equations\paperinfo Springer-Verlag,
1997\endref

\ref\key IT\by J. Inoguchi, M. Toda\paper Timelike Minimal
Surfaces via Loop Groups, Acta Applicandae Mathematicae\paperinfo
Vol.83, Number 3, (2004), 313-355\endref


\ref\key Jo1\by D. Joyce\paper Lectures on Calabi-Yau and special
Lagrangian geometry\paperinfo in:Mark Gross, Daniel Huybrechts and
Dominic Joyce (eds.), Calabi-Yau Manifolds and Related Geometries
(Universitext series), Springer-Berlin, (2003)\endref

\ref\key Jo2\by D. Joyce\paper Special Lagrangian m-folds in $C^m$
with symmetries\paperinfo Duke Mathematical Journal, 115 (2002),
1-51\endref

\ref\key JX\by J. Jost, Y.L. Xin\paper Some aspects of the global
geometry of entire spacelike submanifolds\paperinfo Result. Math.
40 (2001), 233-245\endref

\ref\key Ka\by S. Katayama\paper Global existence for systems of
nonlinear wave equations in two space dimensions, II\paperinfo
Publ. RIMS, Kyoto Univ. 31 (1995), 645-665\endref

\ref\key Kl1\by S. Klainerman\paper Uniform decay estimates and
the Lorentz invariance of the classical wave equation\paperinfo
Comm. Pure Appl. Math.38 (1985), 321-332\endref

\ref\key Kl2\by S. Klainerman\paper The null condition and global
existence to nonlinear wave equations\paperinfo Lectures in
Applied Mathematics 23 (1986), 293-326\endref

\ref\key Ko\by O. Kobayashi\paper Maximal Surfaces in the
3-Dimensional Minkowski Space\paperinfo Tokyo J. Math. 6 (1983),
297C309\endref

\ref\key Li\by H. Lindblad\paper A remark on global existence for
small initial data of the minimal surface equation in Minkowskian
space time\paperinfo Proc. Amer. Math. Soc., 132 (2004),
1095-1102\endref

\ref\key Mi\by T. K. Milnor\paper A conformal analog of
Bernstein's theorem in Minkowski 3-space\paperinfo In: The Legacy
of Sonya Kovalevskaya, Contemp. Math. 64, Amer. Math. Soc.,
Providence, RI, 1987, pp. 123-130\endref

\ref\key Si\by A.C. Silva\paper Lectures on Symplectic
Geometry\paperinfo Lect. Notes in Math., Vol. 1764,
Springer-Verlag, 2001\endref

\ref\key SW1\by R.W. Smyth, T. Weinstein\paper Conformally
homeomorphic Lorentz surfaces need not be conformally
diffeomorphic\paperinfo Proc. Amer. Math. Soc. 123 (1995),
3499-3506\endref

\ref\key SW2\by R.W. Smyth, T. Weinstein\paper How many Lorentz
surfaces are there?\paperinfo In: S. Gindikin (ed.), Topics in
Geometry, in memory of Joseph D'Atri, Birkh\"auser, Basel, 1996,
pp. 315-330\endref

\ref\key Ta\by M.E. Taylor\paper Partial Differential Equations,
Vol.III\paperinfo Applied Mathematical Sciences 117,
Springer-verlag, 1997\endref

\ref\key Wa\by M. Warren\paper Calibrations associated to
Monge-Ampere Equations\paperinfo arXiv:math.AP/0702291,
2007\endref

\endRefs
\enddocument